\documentclass[11pt]{article}
\usepackage{amsthm,amsfonts,amsmath,amssymb}
\newtheorem{thm}{Theorem}[section]

\newtheorem{lem}[thm]{Lemma}
\newtheorem{cor}[thm]{Corollary}

\def\di{\bigm|} \def\lg{\langle} \def\rg{\rangle}

\def\f{\noindent}

\def\mod{\hbox{\rm mod }}

\def\demo{\f{\bf Proof}\hskip10pt}

\def\qed{\hfill $\Box$}

\def\lg{\langle}
\def\rg{\rangle}

\def\rr#1{\item[{\rm (#1)}]}

\def\A{\mathcal{A}$$}

\begin{document}
\title{The Classification of Finite Metahamiltonian $p$-Groups
\thanks{This work was supported by NSFC (No. 11471198 \& 10171006). }}
\author{Xingui Fang\\
School of Mathematics Science, Peking University\\
Beijing, 100087 PR China \\
Lijian An\thanks{Corresponding author. e-mail: anlj@sxnu.edu.cn}\\
Department of Mathematics, Shanxi Normal University\\
Linfen, Shanxi, 041004 PR China}

\date{}
\maketitle

\begin{abstract}
A finite non-abelian group $G$ is called metahamiltonian if every
subgroup of $G$ is either abelian or normal in $G$. If $G$ is
non-nilpotent, then the structure of $G$ has been determined. If
$G$ is nilpotent, then the structure of $G$ is determined by the
structure of its Sylow subgroups. However, the classification of
finite metahamiltonian $p$-groups is an unsolved problem. In this
paper, finite metahamiltonian $p$-groups are completely classified
up to isomorphism.

\end{abstract}

\baselineskip=16pt

\section{Introduction}

To determine a finite group by using its subgroup structure is an
important theme in the group theory. Let $G$ be a finite non-abelian
$p$-group. If every proper subgroup of $G$ is abelian then $G$ is
called {\it minimal non-abelian}, which was classified by Redei \cite{R}.
If every subgroup of $G$ is normal in $G$ then $G$ is called
{\it Hamiltonian}, which was classified by Dedekind \cite{Ded}. The
classifications of minimal non-abelian $p$-groups and Hamiltonian
groups are two classical results in the theory of finite $p$-groups.

As a generalization of minimal non-abelian group, many authors
investigate finite $p$-groups with ¡°many abelian subgroups¡±. Among
these works, the classification of $\mathcal{A}_2$-groups is the
most important one. A finite non-abelian $p$-group $G$ is called an
{\it $\mathcal{A}_2$-group} if $G$ is not minimal non-abelian and all of
its subgroups of index $p$ are either abelian or minimal
non-abelian. Many scholars studied and classified
$\mathcal{A}_2$-groups, see \cite{BJ,Ber2,Dra,Kaz,She,ZSAX}. Resent
years, several important classes of $p$-groups which contain
$\mathcal{A}_2$-group are determined. For example, Xu et al.
\cite{Alj} classified finite $p$-groups all of whose non-abelian
proper subgroups are generated by two elements. An et al.
\cite{AHZ,ALQZ,QYXA,QXA,QZGA} classified finite $p$-groups with a
minimal non-abelian subgroup of index $p$. Zhang et al. \cite{ZZLS}
classified finite $p$-groups all of its subgroups of index $p^3$ are
abelian.

As a generalization of Hamilton groups, many authors investigate
finite $p$-groups with ¡°many normal subgroups¡±. For example,
Passman \cite{Pas} classified finite $p$-groups all of whose
non-normal subgroups are cyclic. Zhang et al. \cite{ZGQX,ZLS,ZS}
classified finite $p$-groups all of whose non-normal subgroups have
orders $\le p^3$.

A non-abelian group $G$ is called {\it metahamiltonian} if every proper
subgroup of $G$ is either abelian or normal in $G$. Obviously,
$\mathcal{A}_2$-groups are metahamiltonian. Groups in
\cite{Pas,ZGQX} are also metahamiltonian. Thus the class of
metahamiltonian $p$-groups is much larger than that of minimal
non-abelian $p$-groups and Hamilton $p$-groups. The classification
of metahamiltonian $p$-groups is an old problem. The present paper
is devoted to the classification.

By the way, Nagrebeckii \cite{Nag2} determined the structure of
finite non-nilpotent metahamiltonian groups. Obviously, a nilpotent
group is metahamiltonian if and only if all its Sylow subgroups are
metahamiltonian. Hence finite metahamiltonian groups are completely
determined.

This paper is divided into four sections. Section 2 is a
preliminary. In section 3, we classify finite metahamiltonian
$p$-groups whose derived group is of exponent $p$, and the case of
exponent $>p$ is dealt  with in section 4.

\smallskip

The sketch of the classification of metahamiltonian $p$-groups is as
follows.

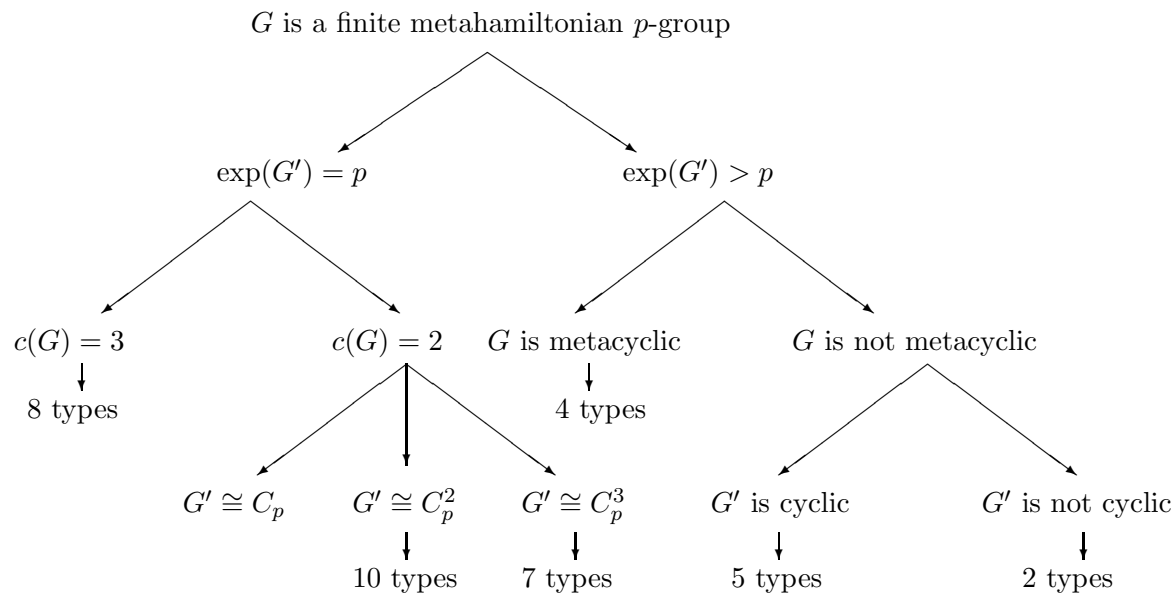
\begin{center}
\vspace{1cm} \setlength{\unitlength}{.90mm}
\begin{picture}(160, 100)(-10, 0)
%\linethickness{1pt} \put(0, 98){\bf Figure 1: $\mathcal{A}_3$-groups
%having an $\mathcal{A}_1$-subgroup of index $p$}

\put(25, 97) {$G$ is a finite metahamiltonian $p$-group}
    \put(60, 94){\vector(-3, -2){22}} \put(60, 94){\vector(3, -2){22}}
    \put(20, 75){$\exp(G')=p$} \put(80, 75){$\exp(G')>p$}
        \put(25, 72){\vector(-4, -3){22}} \put(25, 72){\vector(4, -3){22}}\put(95, 72){\vector(-4, -3){22}} \put(95, 72){\vector(4, -3){22}}
        \put(-10, 50){$c(G)=3$} \put(37, 50){$c(G)=2$}\put(60, 50){$G$ is metacyclic } \put(105, 50){$G$ is not metacyclic }

           \put(48, 48){\vector(4, -3){22}}\put(48, 48){\vector(0, -3){15}}\put(48, 48){\vector(-4, -3){22}}
           \put(15, 26){$G'\cong C_p$}\put(40, 26){$G'\cong C_p^2$} \put(65, 26){$G'\cong C_p^3$}

          \put(48, 23){\vector(0, -3){4}}\put(73, 23){\vector(0, -3){4}}
        \put(40, 15){10 types }\put(65, 15){7 types }

           \put(125, 48){\vector(4, -3){22}}\put(125, 48){\vector(-4, -3){22}}
           \put(93, 26){$G'$ is cyclic} \put(133, 26){$G'$ is not cyclic}

           \put(0, 48){\vector(0, -3){4}}
            \put(-8, 40){8 types }

            \put(75, 48){\vector(0, -3){4}}
            \put(70, 40){4 types }

            \put(103, 23){\vector(0, -3){4}}\put(148, 23){\vector(0, -3){4}}
            \put(96, 15){5 types }\put(139, 15){2 types }

\end{picture}
\end{center}

\section{Preliminaries}

Let $G$ be a finite $p$-group. For a positive integer $t$, $G$ is said to be an {\it
$\A_t$-group} if the greatest index of non-abelian subgroups is $p^{t-1}$.
So ${\cal A}_1$-groups are just the minimal non-abelian $p$-groups.

\medskip

Let $G$ be a finite $p$-group. We define
\begin{center}
 $\Lambda_m(G)=\{ a\in G\di
a^{p^m}=1 \}$, \hskip 8mm $V_m(G)=\{a^{p^m}\di a\in G \}$,

$\Omega_m(G)=\lg
\Lambda_m(G) \rg=\lg a\in G\di a^{p^m}=1\rg$, and
$\mho_m(G)=\lg
V_m(G)\rg=\lg a^{p^m}\di a\in G\rg$.
\end{center}
$G$ is called {\it $p$-abelian}
if $(ab)^{p}=a^{p}b^{p}$ for all $a,b\in G$. We use $c(G)$ and
$d(G)$ to denote the nilpotency class and minimal number of
generators, respectively.

\medskip

We use $C_n$ and  $C_n^m$  to denote the cyclic group and the direct
product of $m$ cyclic groups of order $n$, respectively.
We use $M_p(m,n)$ to denote groups

\begin{center}$\langle a, b \di
a^{p^m}=b^{p^n}=1, a^b=a^{1+p^{m-1}}\rangle$, where $m\geq 2$,
\end{center}
 and
use $M_p(m,n,1)$ to denote groups
\begin{center} $\langle a, b, c \di
a^{p^m}=b^{p^n}=c^p=1, [a, b]=c, [c, a]=[c, b]=1 \rangle$,
 \end{center}
where
 $m+n\geq 3$ for $p=2$ and $m\geq n$.  We can give a presentation of
minimal non-abelian $p$-groups as follows:

\begin{thm}\label{thm=Redei} { \rm (See \cite {R})}%{\rm(R\'edei)}\
\ Let $G$ be a minimal non-abelian $p$-group. Then $G$ is $Q_8$,
$M_p(m,n)$, or $M_p(m,n,1)$.
\end{thm}

A finite $p$-group $G$ is called {\em metacyclic} if it has a cyclic
normal subgroup $N$ such that $G/N$ is also cyclic.

In 1973 King \cite{king} classified metacyclic $p$-groups. In 1988
Newman and Xu (see \cite{NX,XZ1}) found new presentations for these
groups. Theorem~\ref{metacyclic classification} is quoted from \cite{XZ1}.

\begin{thm}
\label{metacyclic classification} $(1)$ Any metacyclic $p$-group
$G$, $p$ odd, has the following presentation:
$$G=\langle a,b\di a^{p^{r+s+u}}=1,\ b^{p^{r+s+t}}=a^{p^{r+s}},\
a^b=a^{1+p^r}\rangle $$ where $r,\ s,\ t,\ u$ are non-negative
integers with $r\ge 1$ and $u\le r$. Different values of the
parameters $r,\ s,\ t$ and $u$ with the above conditions give
non-isomorphic metacyclic $p$-groups. It is denoted to $<r,s,t,u>_p$
in this paper.

$(2)$ Let $G$ be a metacyclic $2$-group. Then $G$ has one of the
following three kinds of presentations:

{\rm(I)} $G$ has a cyclic maximal subgroup. Hence $G$ is dihedral,
semi-dihedral, generalized quaternion, or an ordinary metacyclic
group presented by $$G=\lg a,b\di
a^{2^n}=1,b^2=1,a^b=a^{1+2^{n-1}}\rg.$$

{\rm(II)} Ordinary metacyclic $2$-groups: $$G=\langle a,b\di
a^{2^{r+s+u}}=1,\ b^{2^{r+s+t}}=a^{2^{r+s}}, a^b=a^{1+2^r}\rangle,$$
where $r,$ $s,$ $t,$ $u$ are non-negative integers with $r\ge 2$ and
$u\le r.$ It is denoted to be $<r,s,t,u>_2$ in this paper.

{\rm(III)} Exceptional metacyclic $2$-groups: $$G=\langle a,b\di
a^{2^{r+s+v+t'+u}}=1,\ b^{2^{r+s+t}}=a^{2^{r+s+v+t'}},
a^b=a^{-1+2^{r+v}}\rangle,$$ where $r,$ $s,$ $v,$ $t,$ $t',$ $u$
are non-negative integers with $r\ge 2,$ $t'\le r,$ $u\le 1$,
$tt'=sv=tv=0$, and  if $t'\ge r-1$ then $u=0$. Groups of different
types or of the same type but with different values of parameters
are not isomorphic to each other. It is denoted to be 
$<r,s,v,t,t',u>_2$ in this paper.
\end{thm}

\begin{lem}\label{phi(G')G_3}{\rm (See \cite{Bla})}
\label{metacyclic} Suppose that $G$ is a finite $p$-group. Then $G$
is metacyclic if and only if $G/\Phi(G')G_3$ is metacyclic.
\end{lem}

\begin{lem}\label{metacyclic An}{\rm(See \cite[Lemma J(i)]{Yak0})}
Let $G$ be a metacyclic $p$-group. Then $G$ is an
$\mathcal{A}_n$-group if and only if $|G'|=p^n$.
\end{lem}

In \cite{AZ}, the properties of metahamiltonian $p$-groups are given as follows:

\begin{thm}\label{thm=metahamilton p-gp G, c(G) is at most 3}
Let $G$ be a metahamiltonian $p$-group. Then $c(G)\le 3$. In
particular, $G$ is metabelian.
\end{thm}

\begin{thm}\label{G'<N}Let $G$ be a finite $p$-group. $G$ is metahamiltonian if and only if
$G'$ is contained in every non-abelian subgroup of $G$.
\end{thm}

\begin{thm} \label{thm=d(G)=2--metacyclic}Suppose that $G$ is a finite metahamilton $p$-group. If $d(G)=2$ and $\exp(G')>p$,
then $G$ is metacyclic.
\end{thm}

\begin{thm}\label{thm=c=3 exp(G')=p} Suppose that $G$ is a finite
metahamiltonian $p$-group having an elementary abelian derived
group. If $c(G)=3$, then $G$ is an $\mathcal{A}_2$-group.
\end{thm}

\begin{cor}
\label{cor=c=3 exp(G')=p} Suppose that $G$ is a finite
metahamiltonian $p$-group having an elementary abelian derived
group. If $c(G)=3$, then $d(G)=2$ and $p$ is odd.
\end{cor}

\section{Finite metahamiltonian
$p$-groups whose derived group is of exponent $p$}

In this section, we determine finite
metahamiltonian $p$-groups whose derived group is of exponent $p$. In order to avoid tedious calculations, we provide a proof which relies on some results obtained in other papers. These papers are \cite{ALQZ,AP,QXA,ZSAX}.

\begin{thm}\label{111} Suppose that $G$ is a finite
metahamiltonian $p$-group with $\exp(G')=p$. Then $G$ is one of the following non-isomorphic groups:
\begin{enumerate}

\rr{A} groups with $|G'|=p$.

\rr{B} $c(G)=3$. In this case, $p$ is odd, $d(G)=2$
and $G\in\mathcal{A}_2$.
\begin{enumerate}

\rr{B1} $\lg a_1,b\di
a_1^p=a_2^p=a_3^p=b^{p^m}=1,[a_1,b]=a_{2},[a_2,b]=a_3,[a_3,b]=1,[a_i,a_j]=1\rg,$
where $p\ge 5$ for $m=1$, $p\ge 3$ and $1\leq i,j\leq 3$;

\rr{B2} $\lg a_1,b\di
a_1^p=a_2^p=b^{p^{m+1}}=1,[a_1,b]=a_{2},[a_2,b]=b^{p^m},[a_1,a_2]=1\rg$,
where $p\ge 3$;

\rr{B3} $\lg a_1,b\di
a_1^{p^2}=a_2^p=b^{p^m}=1,[a_1,b]=a_{2},[a_2,b]=a_1^{\nu
p},[a_1,a_2]=1\rg,$ where $p\ge 3$ and $\nu=1$ or a fixed quadratic
non-residue modulo $p$;

\rr{B4} $\lg a_1,a_2,b\di
a_1^9=a_2^3=1,b^3=a_1^3,[a_1,b]=a_2,[a_2,b]=a_1^{-3},[a_2,a_1]=1\rg$.

\rr{B5} $\lg a,b\di a^{p^2}=b^{p^2}=c^p=1,[a,b]=c,[c,a]=b^{\nu p},
[c,b]=a^{p}\rg$, where $p\ge 5$, $\nu$ is a fixed square non-residue
modulo $p$;

\rr{B6} $\lg a,b\di a^{p^2}=b^{p^2}=c^p=1,[a,
b]=c,[c,a]=a^{-p}b^{-lp},[c, b]=a^{-p}\rg$, where $p\ge 5$,
$4l=\rho^{2r+1}-1$, $r=1, 2, \dots, \frac{1}{2}(p-1)$, $\rho$ is the
smallest positive integer which is a primitive root modulo $p$;

\rr{B7} $\lg a, b\di a^9=b^9=c^3=1,[a,b]=c,[c,a]=b^{-3}, [c,
b]=a^3\rg$;

\rr{B8} $\lg a, b\di a^9=b^9=c^3=1, [a,b]=c, [c,a]=b^{-3}, [c,
b]=a^{-3}\rg$.

\end{enumerate}
\rr{C} $c(G)=2$ and $G'\cong C_p^2$.

\begin{enumerate}
 \rr{C1} $K\times A$, where $K=\lg a_1,a_2,b \di a_1^4=a_2^4=1,
b^2=a_1^2,[a_1,a_2]=1, [a_{1},b]=a_{2}^2, [a_{2},b]=a_{1}^2\rg$ and
$A$ is abelian such that $\exp(A)\le 2$;

  \rr{C2} $K\times A$, where $K=\lg a_1,a_2,b,d \di
a_1^4=a_2^4=1, b^2=a_1^2,d^2=a_{2}^2, [a_1,a_2]=1,
[a_{1},b]=a_{2}^2,
[a_{2},b]=a_{1}^2,[a_{1},d]=a_{1}^2,[a_{2},d]=a_{1}^2a_{2}^2,
[b,d]=1\rg$ and $A$ is abelian such that $\exp(A)\le 2$.

  \rr{C3} $K\times A$, where $K=\lg a_1,a_2,a_3\di
 a_1^{p^{m_1+1}}=a_2^{p^{m_2+1}}=a_3^{p^{m_3}}=1,
  [a_1,a_2]=a_1^{p^{m_1}}, [a_1,a_3]=a_2^{p^{m_2}}, [a_2,a_3]=1\rg$,
 $m_1>1$ for $p=2$, $m_1\ge m_2\ge m_3$, and $A$ is abelian such that
$\exp(A)\le p^{m_2}$;

  \rr{C4} $K\times A$, where $K=\lg a_1,a_2,a_3\di
a_1^{p^{m_1+1}}=a_2^{p^{m_2+1}}=a_3^{p^{m_3}}=1,
  [a_1,a_2]=1, [a_1,a_3]=a_2^{p^{m_2}}, [a_2,a_3]=a_1^{\nu p^{m_1}}\rg$,
 $p>2$, $\nu$ is a fixed square non-residue modulo $p$,
$m_1-1=m_2\ge m_3$ or $m_1=m_2\ge m_3$,
 and $A$ is abelian such that
$\exp(A)\le p^{m_2}$;

  \rr{C5} $K\times A$, where $K=\lg a_1,a_2,a_3\di a_1^{p^{m_1+1}}=a_2^{p^{m_2+1}}=a_3^{p^{m_3}}=1,
  [a_1,a_2]=1, [a_1,a_3]=a_2^{p^{m_2}}, [a_2,a_3]=a_1^{kp^{m_1}}a_2^{-p^{m_2}}\rg$,
 $1+4k\not\in (F_p)^2$ for $p>2$, $k=1$ for $p=2$, $m_1=m_2\ge m_3$ and $A$ is abelian such that
$\exp(A)\le p^{m_2}$;

 \rr{C6} $K\times A$, where $K=\lg a_1,a_2,a_3\di
a_1^{p^{m_1+1}}=a_2^{p^{m_2+1}}=a_3^{p^{m_3}}=1,
  [a_1,a_2]=1, [a_1,a_3]=a_2^{p^{m_2}}, [a_2,a_3]=a_1^{p^{m_1}}\rg$,
$m_1-1=m_2\ge m_3$ and $A$ is abelian such that $\exp(A)\le
p^{m_2}$;

 \rr{C7} $K\times A$, where $K=\lg a_1,a_2,a_3\di
a_1^{p^{m_1}}=a_2^{p^{m_2+1}}=a_3^{p^{m_3+1}}=1,
  [a_1,a_2]=a_3^{p^{m_3}}, [a_1,a_3]=a_2^{p^{m_2}}, [a_2,a_3]=1\rg$, $m_1\ge
  m_2=m_3+1$ and $A$ is abelian such that $\exp(A)\le p^{m_3}$;

  \rr{C8} $K\times A$, where $K=\lg a_1,a_2,a_3\di
a_1^{p^{m_1}}=a_2^{ p^{m_2+1}}=a_3^{p^{m_3+1}}=1,
  [a_1,a_2]=a_3^{p^{m_3}}, [a_1,a_3]=a_2^{\nu p^{m_2}}, [a_2,a_3]=1\rg$, $p>2$, $\nu$ is a fixed square non-residue modulo $p$, $m_1\ge
m_2=m_3+1$ or $m_1>m_2=m_3$ and $A$ is abelian such that
$\exp(A)\le p^{m_3}$;

  \rr{C9} $K\times A$, where $K=\lg a_1,a_2,a_3\di a_1^{p^{m_1}}=a_2^{p^{m_2+1}}=a_3^{p^{m_3+1}}=1,
  [a_1,a_2]=a_3^{p^{m_3}}, [a_1,a_3]=a_2^{kp^{m_2}}a_3^{-p^{m_3}}, [a_2,a_3]=1\rg$,
 $1+4k\not\in (F_p)^2$ for $p>2$, $k=1$ for $p=2$, $m_1> m_2=m_3$ and $A$ is abelian such that
$\exp(A)\le p^{m_3}$;

 \rr{C10} $K\times A$, where $K=\lg a_1,a_2,a_3\di
 a_1^{p^{m_1+1}}=a_2^{p^{m_2}}=a_3^{p^{m_3+1}}=1,
  [a_1,a_2]=a_3^{p^{m_3}}, [a_1,a_3]=a_1^{p^{m_1}}, [a_2,a_3]=1\rg$,
$m_1\ge m_2=m_3+1$ and $A$ is abelian such that $\exp(A)\le
p^{m_3}$.

\end{enumerate}

\rr{D} $c(G)=2$ and $G'\cong C_p^3$.
\begin{enumerate}
\rr{D1}  $K\times A$, where $K=\langle a_1, a_2, a_3 \di a_1^{p^{m_1+1}}=a_2^{p^{m_2+1}}=a_3^{p^{m_3+1}}=1,
[a_2,a_3]=a_1^{p^{m_1}},[a_1,a_3]=a_2^{\eta p^{m_2}},[a_1,a_2]=a_3^{ p^{m_3}},[a_3^p,a_1]=[a_3^p,a_2]=1\rangle$,
where $p$ is odd, $m_1=m_2+1=m_3+1$ and $\eta$ is a fixed square
non-residue modulo $p$, and $A$ is abelian with $\exp(A)\le p^{m_3}$;

\rr{D2} $K\times A$, where $K=\langle a_1, a_2, a_3 \di a_1^{p^{m_1+1}}=a_2^{p^{m_2+1}}=a_3^{p^{m_3+1}}=1,
[a_2,a_3]=a_1^{p^{m_1}},[a_1,a_3]=a_2^{l p^{m_2}}a_3^{- p^{m_2}},[a_1,a_2]=a_3^{ p^{m_3}},[a_3^p,a_1]=[a_3^p,a_2]=1\rangle$,
where $p$ is odd, $m_1=m_2+1=m_3+1$ and $1+4l\not\in (F_p)^2$, and $A$ is abelian with $\exp(A)\le p^{m_3}$;

\rr{D3} $K\times A$, where $K=\langle a_1, a_2, a_3 \di a_1^{2^{m_1+1}}=a_2^{2^{m_2+1}}=a_3^{2^{m_3+1}}=1,
[a_2,a_3]=a_1^{2^{m_1}},[a_3,a_1]=a_2^{2^{m_2}},[a_1,a_2]=a_2^{2^{m_2}}a_3^{2^{m_3}},[a_3^2,a_1]=[a_3^2,a_2]=1\rangle$,
where $m_1=m_2+1=m_3+1$, and $A$ is abelian with $\exp(A)\le 2^{m_3}$;

\rr{D4} $K\times A$, where $K=\langle a_1, a_2, a_3 \di a_1^{p^{m_1+1}}=a_2^{p^{m_2+1}}=a_3^{p^{m_3+1}}=1,
[a_2,a_3]=a_1^{p^{m_1}},[a_1,a_3]=a_2^{\eta p^{m_2}},[a_1,a_2]=a_3^{ p^{m_3}},[a_3^p,a_1]=[a_3^p,a_2]=1\rangle$,
where $p$ is odd, $m_1=m_2=m_3+1$ and $\eta$ is a fixed square
non-residue modulo $p$, and $A$ is abelian with $\exp(A)\le p^{m_3}$;

\rr{D5} $K\times A$, where $K=\langle a_1, a_2, a_3 \di a_1^{p^{m_1+1}}=a_2^{p^{m_2+1}}=a_3^{p^{m_3+1}}=1,
[a_2,a_3]=a_1^{p^{m_1}},[a_1,a_3]=a_1^{ p^{m_1}}a_2^{l p^{m_2}},[a_1,a_2]=a_3^{ p^{m_3}},[a_3^p,a_1]=[a_3^p,a_2]=1\rangle$,
where $p$ is odd, $m_1=m_2=m_3+1$ and $1+4l\not\in (F_p)^2$, and $A$ is abelian with $\exp(A)\le p^{m_3}$;

\rr{D6} $K\times A$, where $K=\langle a_1, a_2, a_3 \di a_1^{2^{m_1+1}}=a_2^{2^{m_2+1}}=a_3^{2^{m_3+1}}=1,
[a_2,a_3]=a_1^{2^{m_1}}a_2^{2^{m_2}} ,[a_3,a_1]=a_2^{2^{m_2}},[a_1,a_2]=a_3^{2^{m_3}},[a_3^2,a_1]=[a_3^2,a_2]=1\rangle$,
where $m_1=m_2=m_3+1$, and $A$ is abelian with $\exp(A)\le 2^{m_3}$;

\rr{D7} $K\times A$, where $K=\langle a,b,c \di a^{4}=b^{4}=c^{4}=1,
[b,c]=a^2b^2,[c,a]=b^2c^2,[a,b]=c^2,[c^2,a]=[c^2,b]=1\rangle$, and $A$ is abelian with $\exp(A)\le 2$.
\end{enumerate}
\end{enumerate}
\end{thm}
\demo
By Theorem \ref{thm=metahamilton p-gp G, c(G) is at most 3},
$c(G)\le 3$. If $c(G)=3$, then, by Theorem \ref{thm=c=3 exp(G')=p},
$G\in \mathcal{A}_2$. Checking groups listed in \cite[Lemma 2.4]{AZ}, we
get groups (B1)--(B8). In the following, we may assume that
$c(G)=2$.  Let $N$ be a minimal non-abelian subgroup of $G$. By
Theorem \ref{G'<N}, $G'\le N$. Since $G'\le Z(G)$, $G'\le
\Omega_1(Z(N))=\Omega_1(\Phi(N))$. It follows from Theorem
\ref{thm=Redei} that $G'\le C_p^3$. If $G'\cong C_p$, then $G$ is of
Type (A) in the theorem. If $G'\cong C_p^2$, then, by the following Lemma
\ref{thm=G'=p^2}, $G$ is a group of Type (C1)--(C10) in the theorem. For the case of $G'\cong C_p^3$, Lemma \ref{thm=G'=p^3} gives groups of 
 Type (D1)--(D5) in the theorem. Finally, it is omitted to check that such groups are non-isomorphic metahamiltonian $p$-groups.
\qed

\begin{lem}\label{thm=G'=p^2}Suppose that $G$ is a metahamilton $p$-group. If $G'\cong C_p^2$ and $c(G)=2$,
then $G$ is a group of Type {\rm (C1)--(C10)} as defined in Theorem {\rm \ref{111}}.
\end{lem}
\demo  Let the type of $G/G'$ be $(p^{m_1},p^{m_2},\dots,p^{m_r})$,
where $m_1\ge m_2\ge\dots\ge m_r$. Let
\begin{center}
$G/G'=\lg a_1G'\rg\times\lg
a_2G'\rg\times\dots\times\lg a_rG'\rg$, where $o(a_iG')=p^{m_i}$,
$i=1,2,\dots,r$.
\end{center}
 Then $G=\lg a_1,a_2,\dots,a_r\rg$.

\smallskip

If $m_1=1$, then $G/G'$ is elementary abelian. By Theorem
\ref{G'<N}, $G'\le \lg x,y\rg$ for every non-commutative pair
$x,y\in G$ and hence $\lg x,y\rg$ is minimal non-abelian with order
$p^4$. Such groups were classified in \cite{AP}. By checking the results in \cite{AP},
we get the groups (C3)--(C5) where $m_1=m_2=m_3=1$ and (C1)--(C2). In the
following, we may assume that $m_1>1$.

 Let $i$ be the minimal integer such that
$a_i\not\in Z(G)$. That is, there exists $j>i$ such that
$[a_i,a_j]\neq 1$. If $i\neq 1$, then $a_1\in Z(G)$. Replacing $a_1$
with $a_1a_j$, we get $a_1\not\in Z(G)$. If $i=1$, then we also have
$a_1\not\in Z(G)$.

Let $j$ be the minimal integer such that $[a_1,a_j]\neq 1$. If
$j\neq 2$, then $[a_1,a_2]=1$. Replacing $a_2$ with $a_2a_j$, we get
$[a_1,a_2]\neq 1$. If $j=2$, then we also have $[a_1,a_2]\neq 1$.

Let $k$ be the minimal integer such that $[a_k,a_l]\not\in\lg
[a_1,a_2]\rg$. If $k>2$, then, for all integer $s$, we have
\begin{center}
$[a_1,a_s]\in\lg [a_1,a_2]\rg$ and $[a_2,a_s]\in\lg [a_1,a_2]\rg$.
\end{center}
(1) If $[a_1,a_l]=1$, then, replacing $a_2$ with $a_2a_l$, we have
$[a_2,a_k]\not\in\lg [a_1,a_2]\rg$. (2) If
$[a_1,a_l]=[a_1,a_2]^\alpha$ where $(\alpha,p)=1$, then, letting
$[a_1,a_k]=[a_1,a_2]^\beta$ and replacing $a_2$ with
$a_2a_ka_l^{\alpha^{-1}\beta}$, we have $[a_2,a_l]\not\in\lg
[a_1,a_2]\rg$. Hence we may assume that $k\le 2$.

Let $l$ be the minimal integer such that $[a_k,a_l]\not\in\lg
[a_1,a_2]\rg$. If $l\neq 3$, then $[a_1,a_3]\in\lg [a_1,a_2]\rg$ and
$[a_2,a_3]\in \lg[a_1,a_2]\rg$. Replacing $a_3$ with $a_3a_l$, we
have $[a_k,a_3]\not\in\lg [a_1,a_2]\rg$. Hence we may assume that
$l=3$.

Let $K=\lg a_1,a_2,a_3\rg$. Then $|K'|=|G'|=p^2$. Such groups $K$ were determined in
\cite{ALQZ}. By checking \cite[Table 4]{ALQZ}, $K$ is one of the groups
(C3)--(C10) in Theorem \ref{111}. If $r=3$, then $G=K$. In the following
we may assume that $r\ge 4$.

\medskip

Case 1: $K$ is one of the groups of Type (C3)--(C6) in Theorem \ref{111}.

\medskip

In this case, $G'=\lg a_1^{p^{m_1}},a_2^{p^{m_2}}\rg$ and $[a_1,a_3]=a_2^{p^{m_2}}$. Assume that $a_4^{p^{m_4}}=a_1^{\alpha p^{m_1}}a_2^{\beta p^{m_2}}$. Replacing $a_4$ with $a_4a_1^{-\alpha
p^{m_1-m_4}}$, we have $a_4^{p^{m_4}}=a_2^{\beta p^{m_2}}$ since $m_1>1$.

If $p>2$ or $m_2>1$, then, replacing $a_4$ with $a_4a_2^{-\beta p^{m_2-m_4}}$, we have $a_4^{p^{m_4}}=1$. If
$p=2$ and $m_2=1$, then we claim that there exists
an $x\in\{ a_4,a_4a_2\}$ such that $x^2\in \lg a_1^{2^{m_1}}\rg$.
Otherwise, $a_4^2=a_2^2$. Since $[a_4,a_2]=(a_4a_2)^2\not\in \lg
a_1^{2^{m_1}}\rg$, $\lg a_4,a_2\rg$ is not abelian. It follows from Theorem \ref{G'<N} that
$a_1^{2^{m_1}}\in \lg a_4,a_2\rg$. Hence
$[a_4,a_2]=a_1^{2^{m_1}}a_2^{2}$. Thus $\lg a_4a_2,
a_2a_1^{2^{m_1-1}}\rg$ is neither abelian nor normal in $G$, a
contradiction. Replacing $a_4$ with $x$ or $xa_1^{2^{m_1-1}}$, we
have $a_4^2=1$.

Hence we may assume that $a_4^{p^{m_4}}=1$. We claim that $[a_1,a_4]\in\lg a_2^{p^{m_2}}\rg$. Otherwise, we may
assume that $[a_1,a_4]=a_1^{\gamma p^{m_1}}a_2^{\alpha p^{m_2}}$
where $(\gamma,p)=1$. By calculation, $\lg a_1,a_4a_3^{-\alpha}\rg$
is neither abelian nor normal in $G$, a contradiction. Hence $[a_1,a_4]\in \lg a_2^{p^{m_2}}\rg$.

Let $L=\lg a_1,a_2,a_4\rg$. If $[a_1,a_4]\ne 1$, then, by suitable replacement, we may assume that $[a_1,a_4]=a_2^{p^{m_2}}$. In this case, we claim that $L'=G'$.
 If not, then $L'=\lg a_2^{p^{m_2}}\rg$. Since $G'\not\le \lg a_2,a_4\rg$, $[a_2,a_4]=1$ by Theorem \ref{G'<N}.
Since $K'=G'$, $K'=\lg a_2^{p^{m_2}},[a_2,a_3]\rg$. Hence we may assume that $[a_2,a_3]=a_1^{sp^{m_1}}a_2^{tp^{m_2}}$ where $(s,p)=1$. If $(t,p)=1$, then $\lg a_1^{sp^{m_1-m_2}}a_2^{t}, a_3a_4^{-1}\rg$ is neither abelian nor normal in $G$,
a contradiction. If $t=0$ and $m_1>m_2$, then $\lg a_1a_2, a_3a_4^{-1}\rg$ is neither abelian nor normal in $G$,
a contradiction. If $t=0$ and $m_1=m_2$, then $\lg a_1a_2, a_3a_4^{s-1}\rg$ is neither abelian nor normal in $G$,
also a contradiction.

By a similar argument as above, for $4\le i\le r$, we may assume that
$a_i^{p^{m_i}}=1$ and $[a_1,a_i]=1$ or $a_2^{p^{m_2}}$.
Moreover, we have:
\begin{center}
(*)\ \ If $[a_1,a_i]=a_2^{p^{m_2}}$, then $L'=G'$ where $L=\lg a_1,a_2,a_i\rg$.
\end{center}
For $3\le i<j\le r$, $[a_i,a_j]=1$ by Theorem \ref{G'<N}.

Let $j$ be the maximal integer such that $[a_1,a_j]=a_2^{p^{m_2}}$. Then $[a_1,a_k]=1$ for $j<k\le r$.
For $3\le k<j$, if $[a_1,a_k]=a_2^{p^{m_2}}$, then $[a_1,a_ka_j^{-1}]=1$. Replacing $a_k$ with $a_ka_j^{-1}$ if necessary, we get $[a_1,a_k]=1$.

Let $J=\lg a_1,a_2,a_j\rg$. Then $J$ is one of the groups of Type (C3)--(C6) in Theorem \ref{111} since $J'=\lg a_1^{p^{m_1}},a_2^{p^{m_2}}\rg$.
We claim that $[a_2,a_k]=1$ for $3\le k\le r$ and $k\ne j$. If not, then we will reduce contradictions on two subcases respectively.

Subcase 1: $J$ is the group of Type (C3) in Theorem \ref{111}.

In this subcase, $[a_2,a_j]=1$.
We may assume that $[a_2,a_k]=a_2^{\gamma p^{m_2}}a_1^{\beta p^{m_1}}$ where $(\beta,p)=1$. If $(\gamma,p)=1$, then $\lg a_1^{\beta
p^{m_1-m_2}}a_2^{\gamma},a_k\rg$ is neither abelian nor
normal in $G$, a contradiction. If $\alpha=0$ and $m_1>m_2$,
then $\lg a_1a_2, a_k\rg$ is neither abelian nor normal
in $G$, a contradiction. If $\alpha=0$ and $m_1=m_2$, then $\lg a_1a_2,
a_ka_j^{\beta}\rg$ is neither abelian nor normal in $G$, also
a contradiction.

Subcase 2: $J$ is one of the groups of Type (C4)--(C7) in Theorem \ref{111}.

In this subcase, $[a_1,a_2]=1$ and $G'=\lg [a_2,a_j], a_2^{p^{m_2}}\rg$. Hence we may
assume that $[a_2,a_k]=a_2^{\gamma p^{m_2}}[a_2,a_j]^{\beta}$ where $(\beta,p)=1$.
Let $x=a_k^{-\beta^{-1}}a_j$. Then $[a_1,x]=a_2^{p^{m_2}}$ and $[a_2,x]=a_2^{-\beta^{-1}\gamma p^{m_2}}$.
If $(\gamma,p)=1$, then $\lg a_2,a_ka_j^{-\beta}\rg$ is neither abelian nor
normal in $G$, a contradiction. If $\alpha=0$, then $\lg a_1,a_2,x\rg'=\lg a_2^{p^{m_2}}\rg$. This contradicts (*).

In this case, $G=J\times A$ where $A=\lg a_3\rg\times \dots \times\lg a_{j-1}\rg\times\lg a_{j+1}\rg \times\dots\times\lg a_r\rg$.
Hence we get the groups (C3)--(C7) in Theorem \ref{111}.

\medskip

Case 2: $K$ is one of the groups of Type (C7)--(C10) in Theorem \ref{111}.

\medskip

In this case, $G'=\lg a_s^{p^{m_s}},a_3^{p^{m_3}}\rg$ where $s=1$ or $2$, $[a_1,a_2]=a_3^{p^{m_3}}$ and $[a_2,a_3]=1$.
Assume that $a_4^{p^{m_4}}=a_s^{\alpha p^{m_s}}a_3^{\beta p^{m_3}}$.

If $p>2$ or $m_3>1$, then, replacing $a_4$ with $a_4a_s^{-\alpha p^{m_s-m_4}}a_3^{-\beta p^{m_3-m_4}}$, we have $a_4^{p^{m_4}}=1$. If
$p=2$, $m_3=1$ and $m_s>1$, then, we claim that there exists
an $x\in\{ a_4,a_4a_3\}$ such that $x^2\in \lg a_s^{2^{m_s}}\rg$.
Otherwise, $a_4^2=a_s^{\alpha 2^{m_s}}a_3^2$. Replacing $a_4$ with $a_4a_s^{-\alpha p^{m_s-m_4}}$, we have $a_4^{2}=a_3^{2}$.
Since $[a_4,a_3]=(a_4a_3)^2\not\in \lg
a_s^{2^{m_s}}\rg$, $\lg a_4,a_3\rg$ is non-abelian. It
follows from Theorem \ref{G'<N} that $G'\le \lg a_4,a_3\rg$. Hence
$[a_4,a_3]=a_s^{2^{m_s}}a_3^{2}$. By calculation, $\lg a_4a_3,
a_3a_s^{2^{m_s-1}}\rg$ is neither abelian nor normal in $G$, a
contradiction. Replacing $a_4$ with $x$ or $xa_s^{2^{m_s-1}}$, we
have $a_4^2=1$. If
$p=2$ and $m_s=m_3=1$, then $s=2$ since $m_1>1$. Hence $K$ is a group of Type (C9). In this case, $[a_1,a_3]=a_2^2a_3^2$.
we claim that there exists an
involution in $\{ a_4,a_4a_2,a_4a_3,a_4a_2a_3\}$. Otherwise, since
$a_4^2\ne 1$, we have
\begin{center}
$a_4^2=a_2^2$, $a_3^2$ or $a_2^2a_3^2$.
\end{center}
If $a_4^2=a_3^2$, then, by replacing $a_2,a_3$ with $a_3,a_2a_3$ respectively, it is reduced to $a_4^2=a_2^2$. If $a_4^2=a_2^2a_3^2$, then, by replacing $a_2,a_3$ with $a_2a_3,a_2$ respectively, it is also reduced to $a_4^2=a_2^2$. Hence we may assume that $a_4^2=a_2^2$. Since
$(a_4a_2)^2=[a_4,a_2]\neq 1$, $L=\lg a_4,a_2\rg$ is not abelian. It
follows from Theorem \ref{G'<N} that $G'\le L$. Hence
we may assume that $[a_4,a_2]=a_3^2a_2^{2\alpha}$. If
$[a_4,a_2]=a_3^2a_2^2$, then $\lg a_1a_4,a_2\rg$ is neither abelian
nor normal in $G$, a contradiction. If $[a_4,a_2]=a_3^2$, then
$(a_4a_2)^2=a_3^2$. Since $(a_4a_2a_3)^2\ne 1$,
$[a_4a_2,a_3]=[a_4,a_3]=(a_4a_2a_3)^2\neq 1$. Since $M=\lg
a_4a_2,a_3\rg$ is not abelian, $G'\le \lg a_4a_2,a_3\rg$ by Theorem \ref{G'<N}. Hence we may assume that $[a_4,a_3]=[a_4a_2,a_3]=a_2^2a_3^{2\alpha}$.  Since $(a_4a_3)^2\ne 1$,
$[a_4,a_3]\neq a_2^2a_3^3$. Hence $[a_4,a_3]=a_2^2$. In this case,
$\lg a_1a_4a_2,a_3\rg$ is neither abelian nor normal in $G$, a
contradiction.

By the above argument, we may assume that $a_4^{p^{m_4}}=1$. Let $\{s,t\}=\{ 1,2\}$.  Since $G'\not \le \lg a_t,a_4\rg$, $[a_t,a_4]=1$ by Theorem \ref{G'<N}. By the definition relations of (C7)--(C10), $m_t>m_3$. It follows from Theorem \ref{G'<N} that $[a_ta_3,a_4]=1$ since $G'\not \le \lg a_ta_3,a_4\rg$. Hence $[a_3,a_4]=1$. We claim that $[a_s,a_4]\in\lg a_3^{p^{m_3}}\rg$. Otherwise, we may assume that $[a_s,a_4]=a_s^{\alpha p^{m_s}}a_3^{\beta p^{m_3}}$ where $(\alpha,p)=1$.
By calculation, $\lg a_s,a_t^\beta a_4^{s-t}\rg$ is neither abelian or normal in $G$, a contradiction.
Hence we may assume that $[a_s,a_4]=a_3^{\beta p^{m_3}}$.

We claim that $[a_s,a_4]=1$. If not, then, $(\beta,p)=1$ and we may assume that $[a_s,a_4]=a_3^{p^{m_3}}$ by suitable replacement.
We will reduce contradictions on three subcases respectively.

\medskip

Subcase 1: $s=2,t=1$ and $m_2>m_3$.

In this subcase, $K$ is one of the groups of Type (C7)--(C8).
By the definition relations of Type (C7)--(C8), $[a_1,a_3]=a_2^{\eta p^{m_2}}$ where $\eta=1$ or $\nu$. By calculation,
$\lg a_1 a_4,a_2a_3\rg$ is neither abelian or normal in $G$, a contradiction.

\medskip

Subcase 2: $s=2,t=1$ and $m_2=m_3$.

In this subcase, $K$ is one of the groups of Type (C8)--(C9).
If $K$ is one of the groups of Type (C8), then $[a_1,a_3]=a_2^{\nu p^{m_2}}$. By calculation,
$\lg a_1a_4^{1-\nu},a_2a_3\rg$ is neither abelian or normal in $G$, a contradiction.
If $K$ is one of the groups of Type (C9), then $[a_1,a_3]=a_2^{k p^{m_2}}a_3^{-p^{m_3}}$ where $(k,p)=1$. By calculation,
$\lg a_1 a_4,a_2^ka_3^{-1}\rg$ is neither abelian or normal in $G$, a contradiction.

\medskip

subcase 3: $s=1,t=2$.

In this subcase, $K$ is a group of Type (C10).
By the definition relations of Type (C10), $[a_1,a_3]=a_1^{p^{m_3}}$. By calculation,
$\lg a_1,a_2 a_3a_4^{-1}\rg$ is neither abelian or normal in $G$, also a contradiction.

\medskip

Hence $[a_s,a_4]=1$.
By a similar argument, for $4\le i\le r$, we may assume that $a_i^{p^{m_i}}=1$. Moreover, $[a_1,a_i]=[a_2,a_i]=[a_3,a_i]=1$. For $4\le i<j\le r$, $[a_i,a_j]=1$ by Theorem \ref{G'<N}. In this case, $G=K\times A$ where $A=\lg a_4\rg\times\lg a_5\rg\times \dots\times \lg a_r\rg$.
Hence we get the groups of Type (C7)--(C10) in Theorem \ref{111}.
 \qed

\begin{lem}\label{thm=G'=p^3}Suppose that $G$ is a metahamilton $p$-group. If $G'\cong C_p^3$ and $c(G)=2$,
then $G$ is a group of Type {\rm (D1)--(D7)} as defined in Theorem {\rm \ref{111}}.
\end{lem}

\demo  Let the type of $G/G'$ be $(p^{m_1},p^{m_2},\dots,p^{m_r})$,
where $m_1\ge m_2\ge\dots\ge m_r$, $G/G'=\lg a_1G'\rg\times\lg
a_2G'\rg\times\dots\times\lg a_rG'\rg$, where $o(a_iG')=p^{m_i}$,
$i=1,2,\dots,r$. Then $G=\lg a_1,a_2,\dots,a_r\rg$.
If $[a_i,a_j]\neq 1$, then $G'=\lg a_i^{p^{m_i}},a_j^{p^{m_j}},[a_i,a_j]\rg$ by Theorem \ref{G'<N}.
Hence we have:
\begin{center}
{\bf (*)\ }\ If $a_i^{p^{m_i}}=1$, then $a_i\in Z(G)$.
\end{center}
 Let $i$ be the minimal integer such that
$a_i^{p^{m_i}}\ne 1$. If $i\ne 1$, then $$a_1^{p^{m_1}}=\dots =a_{i-1}^{p^{m_{i-1}}}=1$$ and hence $a_1,\dots a_{i-1}\in Z(G)$ by (*). We claim that $m_i=m_1$. If not, then $(a_1a_j)^{p^{m_1}}=1$ for $j\ge i$. It follows that $a_1a_j\in Z(G)$ by (*) and hence $a_j\in Z(G)$ for $j\ge i$. This contradicts $|G'|=p^3$. Hence we may assume that $a_1^{p^{m_1}}\ne 1$.

Let $j$ be the minimal integer such that $a_j^{p^{m_j}}\not\in \lg a_1^{p^{m_1}}\rg$. If
$j\ne 2$, then we may assume that $a_k^{p^{m_k}}=a_1^{\alpha_k p^{m_1}}$ for $2\le k\le j-1$. By Theorem \ref{G'<N}, $[a_k,a_1]=1$. Replacing $a_k$ with $a_ka_1^{-\alpha_kp^{m_1-m_k}}$, we get $a_k^{p^{m_k}}=1$. By (*), $a_k\in Z(G)$ for $2\le k\le j-1$. We claim that $m_j=m_2$. If not, then $(a_2a_k)^{p^{m_2}}=1$ for $k\ge j$. It follows that $a_2a_k\in Z(G)$ by (*) and hence $a_k\in Z(G)$ for $k\ge j$. This contradicts $|G'|=p^3$. Hence we may assume that $a_2^{p^{m_2}}\not\in \lg a_1^{p^{m_1}}\rg$.

Let $k$ be the minimal integer such that $a_k^{p^{m_k}}\not\in \lg a_1^{p^{m_1}},a_2^{p^{m_2}}\rg$. If $k\ne 3$, then we may assume that $a_w^{p^{m_w}}=a_1^{\alpha_w p^{m_1}}a_2^{\beta_w p^{m_2}}$ for $3\le w\le k-1$.
We claim that $m_k=m_3$. If not, then $m_3>m_k$. Without loss of generality, we may assume that $m_{k-1}>m_k$.
Replacing $a_w$ with $a_wa_1^{-\alpha_w p^{m_1-m_w}}a_2^{\beta_w p^{m_2-m_w}}$, we get $a_w^{p^{m_w}}=1$. By (*), $a_w\in Z(G)$ for $3\le w\le k-1$.
For $w\ge k$, since
$(a_3a_w)^{p^{m_3}}=1$, $a_3a_w\in Z(G)$ by (*). It follows that $a_w\in Z(G)$ for $w\ge k$. This contradicts $|G'|=p^3$. Hence we may assume that $a_3^{p^{m_3}}\not\in \lg a_1^{p^{m_1}},a_2^{p^{p^{m_3}}}\rg$.

If $r=3$, then, by \cite[Theorem 8.1]{QXA}, $G$ is a group of
Type (D1)--(D7) in Theorem \ref{111}. In the following we may
assume that $r\ge 4$.

We claim that there are suitable $a_1,a_2,a_3$ such that the following condition:
\begin{center}
{\bf (**)} For all $x\in G'$, there exists $b\in \lg a_1,a_2,a_3\rg$ such that $x=b^{p^{m_3}}$.
\end{center}
If (**) holds, then for $i>3$, there exists $b_i\in\lg a_1,a_2,a_3\rg$ such that $a_i^{p^{m_i}}=b_i^{p^{m_3}}$. By Theorem \ref{G'<N},
$[a_i,b_i]=1$. Replacing $a_i$ with $a_ib_i^{-p^{m_3-m_i}}$, we get
$a_i^{p^{m_i}}=1$. By (*), $a_i\in Z(G)$. Hence we
get the groups (D1)--(D7) in Theorem \ref{111}.

\medskip

In the following, we prove that we may choose suitable $a_1,a_2,a_3$ satisfying the condition (**).
If $p>2$ or $m_2>1$, then (**) holds. Hence, we only need to deal with the case where $p=2$ and $m_2=1$.

Case 1. $m_1>1$.

If $[a_2,a_3]\ne 1$, then we may assume that $[a_2,a_3]=a_2^{2i}a_3^{2j}a_1^{2^{m_1}}$ by Theorem \ref{G'<N}.
If $[a_2,a_3]=a_2^{2}a_3^{2j}a_1^{2^{m_1}}$, then $\lg a_2a_1^{2^{m_1-1}},a_3\rg$ is neither abelian nor normal in $G$, a contradiction. If $[a_2,a_3]=a_2^{3}a_1^{2^{m_1}}=(a_3a_1^{2^{m_1-1}})^2$, then $\lg a_3a_1^{2^{m_1-1}},a_3\rg$ is neither abelian nor normal in $G$, a contradiction. Hence $[a_2,a_3]=a_1^{2^{m_1}}$. In this case, it is easy to check that $G'=V_1(\lg a_1,a_2,a_3\rg)$.
Hence (**) holds.

Case 2. $m_1=1$.

By an argument similar to the beginning of the proof of Theorem \ref{111},
we may choose suitable $a_1,a_2,a_3$ such that the commutative group of $K=\lg a_1,a_2,a_3\rg$ is of order at least 4.

If there
are two elements in $\{1,a_1,a_2,a_3,a_1a_2,a_1a_3,a_2a_3,a_1a_2a_3\}$ such that the squares are equal to each other, then, by Theorem \ref{G'<N}, they are commutative. It follows that there is an involution in
$\{a_1,a_2,a_3,a_1a_2,a_1a_3,a_2a_3,a_1a_2a_3\}$. By (*), this involution is in the center of $K$, which contradicts $|K'|\ge 4$. Hence $$G'=V_1(K)=\{1,a_1^2,a_2^2,a_3^2,(a_1a_2)^2,(a_1a_3)^2,(a_2a_3)^2,(a_1a_2a_3)^2\}.$$  That is, (**) holds.
\qed

\section{Finite metahamiltonian
$p$-groups whose derived group is of exponent $>p$}

\begin{thm}\label{222}Suppose that $G$ is a finite
metahamiltonian $p$-group with $\exp(G')>p$. Then $G$ is
isomorphic to one of the following non-isomorphic groups:
\begin{enumerate}

\rr{E} $G$ is metacyclic.

\begin{enumerate}
\rr{E1} $\langle a,b\di a^{p^{r+s+u}}=1,\
b^{p^{r+s+t}}=a^{p^{r+s}},\ a^b=a^{1+p^r}\rangle $, where $r\ge 1$,
$u\le r$, $r+1\ge s+u\ge 2$, and if $p=2$ then $r\ge 2$;

\rr{E2} $\lg a,b\di a^{2^{3}}=b^{2^m}=1, a^b=a^{-1}\rg$, where $m\ge
1$;

\rr{E3} $\lg a,b\di a^{2^{3}}=1, b^{2^m}=a^4, a^b=a^{-1}\rg$, where
$m\ge 1$;

\rr{E4} $\lg a,b\di a^{2^{3}}=b^{2^m}=1, a^b=a^3\rg$, where $m\ge
1$.
\end{enumerate}

\rr{F} $G$ is not metacyclic and $G'$ is cyclic and $|G'|\ge p^2$.
\begin{enumerate}
\rr{F1} $K\times A$, where $K=\lg a,b\di
a^{p^{r+s+u}}=1,b^{p^{r+s}}=1,a^b=a^{1+p^{r}}\rg$, $u\le r$, $r+1>
s+u\ge 2$, and $A\ne 1$ is abelian such that $\exp(A)\le p^{
(r+1)-(s+u)}$;

\rr{F2} $K\times A$, where $K=\lg a,b\di
a^{p^{r+t+u}}=1,b^{p^{r}}=1,a^b=a^{1+p^{r+t}}\rg$, $t\ge 1$, $r\ge
u\ge 2$, and $A\ne 1$ is abelian such that $\exp(A)\le
p^{t+(r+1)-u}$;

\rr{F3} $K\times A$, where $K=\lg a,b\di
a^{p^{r+s}}=1,b^{p^{r+s+t}}=1,a^b=a^{1+p^{r}}\rg$, $t\ge 1$, $r+1>
s\ge 2$, and $A\ne 1$ is abelian such that $\exp(A)\le p^{(r+1)-s}$;

\rr{F4} $K\times A$, where $K=\lg a,b\di
a^{p^{r+s+u}}=1,b^{p^{r+s+t}}=a^{p^{r+s}},a^b=a^{1+p^{r}}\rg$,
$stu\neq 0$, $r+1> s+u\ge 2$, and $A\ne 1$ is abelian such that
$\exp(A)\le p^{(r+1)-(s+u)}$;

\rr{F5} $(K\rtimes B)\times A$, where $K=\lg a,b\di
a^{p^{r+t+u}}=1,b^{p^{r}}=1,a^b=a^{1+p^{r+t}}\rg$, $B=\lg
b_1\rg\times\lg b_2\rg\times\dots\times\lg b_f\rg$ such that
$o(b_i)=p^{r_i}$, $[a,b_i]=a^{p^{r+t_i}}$, $[b,b_i]=1$,
$\max\{t,u-2\}<t_1<t_2<\dots<t_f<t+u$,
$r+t>r_1+t_1>r_2+t_2>\dots>r_f+t_f\ge t+u\ge t+2$, and $A$ is
abelian such that $\exp(A)\le p^{t+(r+1)-u}$.
\end{enumerate}
\rr{G} the type of $G'$ is $(p^\alpha,p)$ where $\alpha\ge 2$.

\begin{enumerate}
\rr{G1} $\lg a_1,a_2,a_3\di
 a_1^{p^{m_1+1+m_2}}=a_2^{p^{m_2+1}}=a_3^p=1,\
  [a_1,a_2]=a_1^{p^{m_1}},\  [a_1,a_3]=a_2^{p^{m_2}}, [a_2,a_3]=1\rg$,
  where $p>2$ and $m_1>
  m_2\ge 1$;

\rr{G2} $K\times A$, where $K=\lg a_1,a_2,a_3\di
 a_1^{p^{m_1+1+k}}=a_2^{p^{m_2+1}}=a_3^{p^{m_3}}=1,
  [a_1,a_2]=a_1^{p^{m_1}}, [a_1,a_3]=a_2^{p^{m_2}}, [a_2,a_3]=1\rg$, $m_1\ge m_2\ge m_3$, $1\le
  k\le\min\{m_1-m_3,m_2-m_3+1,m_2-1\}$ and $A$ is abelian such that $\exp(A)\le p^{m_2-k}$.
\end{enumerate}
\end{enumerate}
\end{thm}
\demo If $G$ is metacyclic, then, by Lemma \ref{metacyclic
metahamilton}, $G$ is a group of Type (E1)--(E4) in the theorem. In the
following, we may assume that $G$ is not metacyclic. If $G'$ is
cyclic, then, by Lemma \ref{thm=G'=C_p^n}, $G$ is a group of Type
(F1)--(F5) in the theorem. If $G'$ is not cyclic, then, by Lemma
\ref{thm=G'=C_p^n*c_p}, $G$ is a group of Type (G1)--(G2) in
theorem. Finally, it is omitted to check that such groups are non-isomorphic metahamiltonian $p$-groups.\qed

\begin{lem}\label{metacyclic metahamilton}Suppose that $G$ is a metacyclic $p$-group and $|G'|\ge p^2$.
If $G$ is metahamiltonian, then $G$ is a group of Type
{\rm (E1)--(E4)} as defined in Theorem {\rm \ref{222}}.
\end{lem}

\demo Case 1: $p>2$ or $G$ is an ordinary metacyclic $2$-group. That
is,
\begin{center}$G=\lg a,b\di
a^{p^{r+s+u}}=1,b^{p^{r+s+t}}=a^{p^{r+s}},a^b=a^{1+p^r}\rg$,\end{center} where
$r\ge1$, $u\le r$, and if $p=2$ then $r\ge 2$.

Since $|G'|\ge p^2$, we have $s+u\ge 2$. We only need to prove that
$r+1\ge s+u$. Otherwise, $r+1<s+u$. By calculation,\begin{center} $[a^{p^{r+1}},
b]=a^{-p^{r+1}}(a^{p^{r+1})})^b=a^{p^{2r+1}}\neq 1$. \end{center}Hence
$\lg a^{p^{r+1}}, b\rg$ is neither abelian nor normal in $G$, a contradiction.
Thus $r+1\ge s+u$ and $G$ is a group of Type (E1) in Theorem \ref{222}.

\medskip

Case 2: $p=2$ and $G$ is not an ordinary metacyclic $2$-group.

Let $o(a)=2^n$ and $H=\lg a^{2^{n-2}},b\rg$. Since $H'=\lg
a^{2^{n-1}}\rg$, $H$ is not abelian. It follows that $H\unlhd G$. By
Theorem \ref{G'<N}, $a^2\in H$. Hence $n\le 3$ and $|G'|=4$. By
Lemma \ref{metacyclic An}, $G\in\mathcal{A}_2$. By \cite[Lemma 2.4]{AZ},
we get groups of (E2)--(E4) in Theorem \ref{222}. \qed

\medskip

We need the following two lemmas on number theory. Proofs are omitted.

\begin{lem}\label{lem=num p odd}
Suppose that $U=U(p^{n})$ is the multiplicative group containing of
all the invertible elements of $\Bbb Z/p^{n}\Bbb Z$, where $p$ is an odd
prime and $n$ is a positive integer. That is,
$$U=\{x\in \Bbb Z/p^{n}\Bbb Z\mid(x,p)=1\}.$$
 Let $S(U)\in {\rm Syl}_p(U)$. Then
 $$S(U)=\{x\in U\mid x\equiv1\ (\mod p)\},$$
and $S(U)$ is cyclic with order $p^{n-1}$. $S_{i}(U)$ where $0\leq
i<n$, the unique subgroup of $S(U)$ of order $p^{i}$, is
 $$S_{i}(U)=\{x\in U\mid x\equiv1\ (\mod p^{n-i})\}.$$
\end{lem}

\begin{lem}\label{lem=num p 2}
Suppose that $U=U(2^{n})$ is the multiplicative group containing of
all invertible elements of $\Bbb Z/2^{n}\Bbb Z$, where $n\ge 2$ is a
positive integer. Then
$$\begin{array}{cl}
  U & =\lg -1\rg\times\lg 1+2^2\rg(\cong C_2\times C_{2^{n-2}}) \\
   & =\{ \varepsilon +i2^m\di \varepsilon=\pm 1, 2\le m\le n, 1\le i\le 2^{n-m}\ {\rm and }\ i\ {\rm is \ odd}\}
\end{array}$$
For $m<n$, the order of $\varepsilon+i2^m$ is $2^{n-m}$ and $\lg
\varepsilon+i2^m\rg=\lg \varepsilon+j2^m\rg$ for all odd $j$.
 \end{lem}

\begin{lem}\label{thm=G'=C_p^n}Suppose that $G$ is a metahamilton $p$-group and $G$ is not metacyclic.
If $|G'|\ge p^2$ and $G'$ is cyclic, then $G$ is a group of Type
{\rm (F1)--(F5)} in Theorem {\rm \ref{222}}.
\end{lem}
\demo By Theorem \ref{thm=d(G)=2--metacyclic}, $d(G)>2$. Let $G'=\lg
c\rg$, the type of $G/G'$ be $(p^{m_1},p^{m_2},\dots,p^{m_w})$ where
$m_1\ge m_2\ge\dots\ge m_w$.
Let
\begin{center}$G/G'=\lg a_1G'\rg\times\lg
a_2G'\rg\times\dots\times\lg a_wG'\rg$ where $o(a_iG')=p^{m_i}$,
$i=1,2,\dots,w$.\end{center} Then $G=\lg a_1,a_2,\dots,a_w\rg$.

Let $i$ be the minimal integer such that $a_i\not\in
C_G(G/\mho_1(G'))$. Then there exists $j>i$ such that $G'=\lg
[a_i,a_j]\rg$. If $i\neq 1$, then $a_1\in C_G(G/\mho_1(G'))$.
Replacing $a_1$ with $a_1a_j$, we have $G'=\lg [a_1,a_i]\rg$.

Let $j$ be the minimal integer such that $G'=\lg [a_1,a_j]\rg$. If
$j\neq 2$, then $[a_1,a_2]\in\mho_1(G')$. Replacing $a_2$ with
$a_2a_j$, we have $G'=\lg [a_1,a_2]\rg$.

Let $K=\lg a_1,a_2\rg$. By Theorem \ref{thm=d(G)=2--metacyclic}, $K$
is metacyclic. Hence $K$ is one of the groups in Theorem
\ref{metacyclic metahamilton}. That is, $K$ is one of the groups
(E1)--(E4) in Theorem \ref{222}.

\medskip

Step 1: We claim that $K$ is one of the groups of Type (E1) in Theorem \ref{222}.

\medskip

If not, then we may assume that $K=\lg a,b\rg$ satisfying the
relations of Type (E2)--(E4) in Theorem \ref{222}. That is,
\begin{center}
$a^{2^3}=1,b^{2^m}\in \mho_1(K')=\lg a^4\rg$ and $[a,b]\equiv a^2\ (\mod \mho_1(K'))$.\end{center}
Obviously, $G'=K'=\lg a^2\rg$ and $m_3=m_4=\dots=m_w=1$.

\medskip

Case 1: $a_3^2\in \mho_1(K')$ and $[a_3,b]\in \mho_1(K')$.

If $[a_3,b]=a^4$, then $\lg a_3,b\rg$ is neither abelian nor normal
in $G$, a contradiction. If $[a_3,b]=1$, then $\lg a_3a^2,b\rg$ is
neither abelian nor normal in $G$, a contradiction.

Case 2: $a_3^2\in \mho_1(K')$ and $[a_3,b]\equiv a^2 \ (\mod \mho_1(K'))$.

If $[a_3,a]\equiv a^2 \ (\mod \mho_1(K'))$, then $(a_3a)^2\in \mho_1(K')$ and
$[a_3a,b]\in \mho_1(K')$. Replacing $a_3$ with $a_3a$, it is reduced to
Case 1. Hence we may assume that $[a_3,a]\in \mho_1(K')$. Since
$[a_3,a^2]=[a_3,a]^2=1$,  $[a_3,G']=1$. By calculation,
$1=[a_3^2,b]=[a_3,b]^2[a_3,b,a_3]=[a_3,b]^2$. Hence $[a_3,b]\in
\mho_1(K')$, a contradiction.

Case 3: $a_3^2\equiv a^2 \ (\mod \mho_1(K'))$.

If $[a_3,a]\in \mho_1(K')$, then, replacing $a_3$ with $a_3a$, it is
reduced to Case 1 or Case 2. Hence we may assume that $[a_3,a]\equiv
a^2\ (\mod \mho_1(K'))$. Since $a_3^2\equiv a^2\ (\mod \mho_1(K'))$,
$[a_3^2,b]=[a^2,b]=a^4$. It follows that $[a_3,b]\equiv a^2\  (\mod
\mho_1(K'))$. Since $(a_3a)^2\equiv a^2\  (\mod \mho_1(K'))$, similar reason as
above gives that $[a_3a,b]\equiv a^2\  (\mod \mho_1(K'))$. Hence $[a,b]\in
\mho_1(K')$, a contradiction.

\medskip

Step 2: By suitable replacement, we may assume $a_i^{p^{m_i}}=1$, where
$3\le i\le w$. Moreover, $[a_i,a_j]=1$ for all $3\le i,j\le w$.

\medskip

By Step 1, $K\cong <r,s,t,u>_p$ where $r\ge 1$, $u\le r$, $r+1\ge
s+u$, and if $p=2$ then $r\ge 2$. Assume that $$K=\lg a,b\di
a^{p^{r+s+u}}=1,b^{p^{r+s+t}}=a^{p^{r+s}},a^b=a^{1+p^r}\rg.$$
Let $L=\lg a,a_i\rg$ and $x_i\in L$ such that $L=\lg
a,x_i\rg$ and $\lg x_i\rg\cap\lg a\rg$ has minimal order. We claim
that $x_i^{p^{m_i}}=1$. Otherwise, we may assume that
\begin{center}$\lg
x_i\rg\cap \lg a\rg=\lg a^{p^\alpha}\rg$ and $\lg [x_i,a]\rg=\lg
a^{p^\beta}\rg$ where $\alpha\ge r$ and $\beta\ge r$.\end{center}
 Then there
exist integers $y$ and $z$ such that $(yz,p)=1$,
$x_i^{p^{m_i}}=a^{yp^{\alpha}}$ and $[x_i,a]=a^{zp^{\beta}}$. By
calculation,
\begin{eqnarray*}
(x_ia^{-yp^{\alpha-m_i}})^{p^{m_i}}&=&x_i^{p^{m_i}}[x_i,a^{yp^{\alpha-m_i}}]^{p^{m_i}\choose
2}[x_i,a^{yp^{\alpha-m_i}},x_i]^{p^{m_i}\choose
3}a^{-yp^{\alpha}}\\
&=&a^{yzp^{\alpha+\beta-m_i}{p^{m_i}\choose
2}}[a^{yzp^{\alpha+\beta-m_i}{p^{m_i}\choose 3}},x_i]\end{eqnarray*}
 Noting that
$\beta\ge r\ge 2$ for $p=2$, we have
$(x_ia^{-yp^{\alpha-m_i}})^{p^{m_i}}\in\lg a^{p^{\alpha+1}}\rg$,
which is contrary to the choice of $x_i$. Replacing $a_i$ with
$x_i$,
we have $a_i^{p^{m_i}}=1$ where $3\le i\le w$.

For $3\le
i,j\le w$, we claim that $[a_i,a_j]=1$. Otherwise,
Theorem \ref{G'<N} gives that $G'\le \lg a_i,a_j\rg$. It is easy to see that $\lg a_i,a_j\rg$ is not metacyclic. This contradicts Theorem \ref{thm=d(G)=2--metacyclic}.

\medskip

Step 3: $K$ is one of the following groups:
\begin{enumerate}
\rr{A} $\lg a,b\di
a^{p^{r+s+u}}=1,b^{p^{r+s}}=1,a^b=a^{1+p^{r}}\rg$, where $r\ge 2$
for $p=2$ and $r+1\ge s+u\ge 2$;

\rr{B} $\lg a,b\di
a^{p^{r+t+u}}=1,b^{p^{r}}=1,a^b=a^{1+p^{r+t}}\rg$, where $t\ge 1$
and $r\ge u\ge 2$;

\rr{C} $\lg a,b\di
a^{p^{r+s}}=1,b^{p^{r+s+t}}=1,a^b=a^{1+p^{r}}\rg$, where $r\ge 2$
for $p=2$, $t\ge 1$ and $r+1\ge s\ge 2$;

\rr{D} $\lg a,b\di
a^{p^{r+s+u}}=1,b^{p^{r+s+t}}=a^{p^{r+s}},a^b=a^{1+p^{r}}\rg$, where
$r\ge 2$, $stu\neq 0$ and $r+1\ge s+u\ge 2$.
\end{enumerate}

\medskip

Assume that $K=\lg a,b\di
a^{p^{r+s+u}}=1,b^{p^{r+s+t}}=a^{p^{r+s}},a^b=a^{1+p^r}\rg$. If
$t=0$, then we have $(ba^{-1})^{p^{r+s}}=1$ for $p>2$ and
$(ba^{2^u-2^{r-1}-1})^{2^{r+s}}=1$ for $p=2$. Replacing $b$ with
$ba^{-1}$ or $ba^{2^u-2^{r-1}-1}$ respectively, we get a group of
Type (A). In the following we may assume that $t\ge 1$. If $s=0$,
then $(a^{-1}b^{p^t})^{p^r}=1$. Replacing $a$ and $b$ with $b$ and
$a^{-1}b^{p^t}$, respectively, we get a group of Type (B). If $u=0$,
then we get a group of Type (C). If $su\neq 0$, then we get a group
of Type (D).

\medskip

Step 4: Determine $G$ in which $K$ is a direct factor. That is,
$G=K\times A$. Since $K'=G'$, $A$ is abelian.

Case 1: $K$ is a group of Type (A) in Step 3.

Let $d\in A$ and $o(d)=p^e$. By calculation,
$$[a^{p^{s+u-1}}d,b]=a^{p^{r+s+u-1}}\neq 1.$$ It follows that
$$a^{p^r}\in \lg (a^{p^{s+u-1}}d)^{p^e}\rg=\lg a^{p^{e+s+u-1}}\rg.$$
Hence $e+s+u-1\le r$. By the arbitrariness of $d$, we get
$\exp(A)\le p^{(r+1)-(s+u)}$. Since $G$ is not metacyclic,
$A\ne 1$. It follows that $r+1>s+u$. Hence we get a group of Type
(F1) in Theorem \ref{222}.

Case 2: $K$ is a group of Type (B) in Step 3.

Let $d\in A$ and $o(d)=p^e$. By calculation,
$$[a^{p^{u-1}}d,b]=a^{p^{r+t+u-1}}\neq 1.$$ It follows that
$$a^{p^{r+t}}\in \lg (a^{p^{u-1}}d)^{p^e}\rg=\lg a^{p^{e+u-1}}\rg.$$
Hence $e+u-1\le r+t$. By the arbitrariness of $d$, we get
$\exp(A)\le p^{t+(r+1)-u}$. Hence $G$ is a group of Type (F2) in Theorem \ref{222}.

Case 3: $K$ is a group of Type (C) or (D) in Step 3.

Let $d\in A$ and $o(d)=p^e$. By calculation,
$$[a^{p^{s+u-1}}d,b]=a^{p^{r+s+u-1}}\neq 1.$$ It follows that
$$a^{p^r}\in \lg (a^{p^{s+u-1}}d)^{p^e}\rg=\lg a^{p^{e+s+u-1}}\rg.$$
Hence $e+s+u-1\le r$. By that arbitrariness of $d$, we
get $\exp(A)\le p^{(r+1)-(s+u)}$. Since $G$ is not metacyclic,  $A\ne 1$. It follows that $r+1>s+u$. Hence we get a group of
Type (F3) or (F4) in Theorem \ref{222}.

\medskip

Step 5: Determine $G$ in which $K$ is not a direct factor.

\medskip

Let $G=H\times A$, where $K<H$ and $A$ is as large as possible for $K$. Since $K'=G'$, $A$ is abelian. By Step
2, we may assume that $H=K\rtimes B$ where $B=\lg b_1\rg\times\lg
b_2\rg\times\dots\times\lg b_f\rg$ such that $o(b_i)=p^{r_i}$, $o(bG')\ge
r_1\ge r_2\ge\dots\ge r_f$.

We claim that $K$ is neither a group of Type (C) nor
(D) in Step 3. Otherwise, by calculation, $\lg ab^{-p^t}\rg\cap \lg a\rg=1$.
Since $G'\not\le \lg ab^{-p^t},b_i\rg$, Theorem \ref{G'<N} gives
that $[ab^{-p^t},b_i]=1$. Similar reason gives that $[b,b_i]=1$.
Hence $H=K\times B$, which is contrary to the choice of $H$.

If $K$ is a group of Type (A) in Step 3, then we claim that $s=0$. Otherwise,
by calculation, $\lg ab \rg\cap\lg a\rg\le\lg a^{p^{r+1}}\rg$. Since
$G'\not\le \lg ab,b_i\rg$, Theorem \ref{G'<N} gives that
$[ab,b_i]=1$. Similar reason gives that $[b,b_i]=1$. Hence
$H=K\times B$, which is contrary to the choice of $H$.

By the above argument, we may assume that $$K=\lg a,b\di
a^{p^{r+t+u}}=1,b^{p^{r}}=1,a^b=a^{1+p^{r+t}}\rg,$$ where $t\ge 0$
and $r\ge u\ge 2$. Since $G'\not\le \lg b,b_i\rg$, Theorem
\ref{G'<N} gives that $[b,b_i]=1$.

Let $j$ be the minimal positive integer such that $[a,b_i]$ has
maximal order. We may assume that $j=1$, replacing $b_1$ with
$b_1b_j$ when it is necessary. Similarly, we may assume that $\lg
[a,b_1]\rg\ge\lg [a,b_2]\rg\ge\dots\ge\lg [a,b_f]\rg$.

Assume that $[a,b_i]=a^{\gamma_ip^{r+t_i}}$ where $(\gamma_i,p)=1$.
Then $t\le t_1\le t_2\le\dots\le t_f$. Note that
$a^b=a^{1+\gamma_ip^{r+t_i}}$. By Lemma \ref{lem=num p odd} and
\ref{lem=num p 2}, there exists positive integer $w$ such that
$$(1+\gamma_ip^{r+t_i})^j\equiv 1+p^{r+t_i} \ (\mod p^{r+t+u}).$$
Replacing $b_i$ with $b_i^w$, we have $[a,b_i]=a^{p^{r+t_i}}$.

Case 1: $t_1>t$.

If $t_2=t_1$, then $b_1b_2^{-1}$ is a direct factor of $H$, a
contradiction. So $t_1<t_2$. Similarly, we have $$t<t_1<t_2<\dots<t_f.$$ If
$(b_1b^{-p^{t_1-t}})^{p^{r_1}}=1$, then $b_1b^{-p^{t_1-t}}$ is a
direct factor of $H$, a contradiction. Hence
$(b_1b^{-p^{t_1-t}})^{p^{r_1}}\neq 1$. It follows that
$b^{p^{r_1+t_1-t}}\neq 1$. Hence $r_1+t_1-t<r$. Thus
$$r-r_1>t_1-t>0.$$ Similarly, we have $$r_i+t_i>r_{i+1}+t_{i+1}.$$ By
Lemma \ref{lem=num p odd} and \ref{lem=num p 2}, in the
multiplicative group consisting of all invertible elements of $\Bbb
Z/p^{r+t+u}\Bbb Z$, the order of $1+p^{r+t_f}$ is $p^{t+u-t_f}$.
Since $[a,b_f^{p^{r_f}}]=1$, we have
$$a^{b_f^{p^{r_f}}}=a^{(1+p^{r+t_f})^{p^{r_f}}}=a.$$ It follows that
$r_f\ge t+u-t_f$. Thus $$t_f+r_f\ge t+u.$$
By calculation, $$\lg ba^{p^{t-t_1+u-1}}\rg\cap\lg a\rg=\lg
(ba^{p^{t-t_1+u-1}})^{p^r}\rg=\lg a^{p^{r+t-t_1+u-1}}\rg.$$ Let
$N=\lg ba^{p^{t-t_1+u-1}},b_1\rg$. Since
$[ba^{p^{t-t_1+u-1}},b_1]=a^{p^{r+t+u-1}}\neq 1$, $N$ is not
abelian. By Theorem \ref{G'<N}, $G'\le N$. It follows that
$r+t-t_1+u-1\le r+t$. Thus $$t_1\ge u-1.$$
Finally, by an argument similar to Step 4, we have $\exp(A)\le
p^{t+(r+1)-u}$. Hence we get a group of Type (F5) in Theorem \ref{222}.
In this case, $\exp(A)\le p^r$.

Case 2: $t_1=t$.

Suppose that $h$ is the maximal positive integer such that $t_h=t$. Let $r'=r_h$, $t'=t+(r-r_h)$ and
$\tilde{K}=\lg a,b_h\rg$. Then $$\tilde{K}=\lg a,b_h\di
a^{p^{r'+t'+u}}=1,b_h^{p^{r'}}=1,a^{b_h}=a^{1+p^{r'+t'}}\rg.$$
If $h<f$, then we let $f'=f-h$. For $1\le
i\le f'$, let
\begin{center}
$b_i'=b_{h+i}$, $t_i'=t_{h+i}$,
$\tilde{B}=\lg b_1'\rg\times\dots\lg b_{f'}'\rg$, $\tilde{H}=\tilde{K}\rtimes\tilde{B}$, and

$\tilde{A}=A\times \lg bb_h^{-1}\rg\times\lg b_1b_h^{-1}\rg\dots\lg  b_{h-1}b_h^{-1}\rg$.
\end{center}
Then $G=\tilde{H}\times \tilde{A}$, where $\tilde{A}$ is as large as possible.
Notice that $t_1'>t'$. By a similar argument to Case 1, we get
a group of Type (F5) in Theorem \ref{222}.

If $h=f$, then we also have $G=\tilde{H}\times \tilde{A}$. The difference in this case from the case $h<f$ is $\tilde{H}=\tilde{K}$. By an argument similar to Step 4, we
have $\exp(A)\le p^{t'+(r'+1)-u}$. Hence we get a group of Type (F2)
in Theorem \ref{222}.
\qed

\begin{lem}\label{thm=G'=C_p^n*c_p} Suppose that $G$ is a finite metahamilton $p$-group. If $\exp(G')>p$ and $G'$ is not cyclic,
then $G$ is a group of Type $(G1)$--$(G2)$ in Theorem {\rm \ref{222}}.
\end{lem}
\demo Let $H\le G$ such that $d(H)=2$ and $\exp(H')>p$. By Theorem
\ref{thm=d(G)=2--metacyclic}, $H$ is metacyclic. By Theorem
\ref{G'<N}, $G'<H$ and hence $G'$ is metacyclic.

Let $N=\mho_1(G')$ and $\bar{G}=G/N$. Then $\bar{G}'\cong C_p^2$. By
Theorem \ref{thm=d(G)=2--metacyclic}, $d(G)>2$ and hence
$d(\bar{G})>2$. By Corollary \ref{cor=c=3 exp(G')=p},
$c(\bar{G})=2$. Hence $\bar{G}$ is a group in Theorem
\ref{thm=G'=p^2}. That is, $\bar{G}$ is a group of Type (C1)--(C10)
in Theorem \ref{111}.

Suppose that $\bar{G}$ is a group of Type (C1) in Theorem \ref{111}. That
is, $\bar{G}=\bar{K}\times \bar{A}$, where $$\bar{K}=\lg
\bar{a}_1,\bar{a}_2,\bar{b} \di \bar{a}_1^4=\bar{a}_2^4=1,
\bar{b}^2=\bar{a}_1^2,[\bar{a}_1,\bar{a}_2]=1,
[\bar{a}_{1},\bar{b}]=\bar{a}_{2}^2,
[\bar{a}_{2},\bar{b}]=\bar{a}_{1}^2\rg$$ and
  $\bar{A}$ is abelian such that $\exp(\bar{A})\le
  2$. Then
  \begin{center}$G'=\lg [a_1,b],[a_2,b],
  \mho_1(G')\rg=\lg a_1^2,a_2^2\rg$ and $\mho_1(G')=\lg
  a_1^4,a_2^4\rg$.\end{center}
Let $M$ be a maximal subgroup of $\mho_1(G')$ such that $M\unlhd G$.
Then we may assume that
\begin{center}$M=\lg e,
  \mho_2(G')\rg$, $[a_1,a_2]\equiv e^i \
(\mod M)$,

$b^2\equiv a_1^2e^j \ (\mod M)$ and $[a_1,b]\equiv
a_2^2e^k \ (\mod M)$. \end{center}
It follows from $[a_1,a_2]\equiv e^i \ (\mod
M)$ that $[a_1^2,a_2]\equiv [a_1,a_2^2]\equiv 1 \ (\mod M)$. It
follows from $b^2\equiv a_1^2e^j \ (\mod M)$ that $[a_1^2,b]\equiv
[a_1,b^2]\equiv 1 \ (\mod M)$. On the other hand, it follows from
$[a_1,b]\equiv a_2^2e^k \ (\mod M)$ that $[a_1^2,b]\equiv
[a_1,b]^2[a_1,b,a_1]\equiv a_2^4 \ (\mod M)$. It follows that
$a_2^4\in M$ and hence $M=\lg a_1^8, a_2^4\rg$.

Let $L=\lg a_1M,bM\rg$. Since $\exp(L')=2$, Theorem \ref{thm=c=3
exp(G')=p} gives that $c(L)=2$. It follows that $[a_2^2,b]\equiv 1 \
(\mod M)$. On the other hand, $[a_2^2,b]\equiv
[a_2,b]^2[a_2,b,a_2]\equiv a_1^4 \ (\mod M)$. It follows that
$a_1^4\in M$. Hence $M=\mho_1(G)$, a contradiction.

Similar reasoning gives that $\bar{G}$ is not a group of Type (C2) in Theorem \ref{111}.

Suppose that $\bar{G}$ is a group of Type (C4) in Theorem \ref{111}. That
is, $\bar{G}=\bar{K}\times \bar{A}$, where
\begin{center}
$\bar{K}=\lg
\bar{a}_1,\bar{a}_2,\bar{a}_3\di
\bar{a}_1^{p^{m_1+1}}=\bar{a}_2^{p^{m_2+1}}=\bar{a}_3^{p^{m_3}}=1,
  [\bar{a}_1,\bar{a}_2]=1, [\bar{a}_1,\bar{a}_3]=\bar{a}_2^{p^{m_2}},
 $

 $ [\bar{a}_2,\bar{a}_3]=\bar{a}_1^{\nu p^{m_1}}\rg$, $p>2$, $\nu$ is a fixed square non-residue modulo $p$,

$m_1-1=m_2\ge m_3$ or $m_1=m_2\ge m_3$,
 and
 $\bar{A}$ is abelian such that
$\exp(\bar{A})\le p^{m_2}$.
\end{center}
Then $G'=\lg [a_1,a_3],[a_2,a_3],
  \mho_1(G')\rg=\lg [a_1,a_3],[a_2,a_3]\rg=\lg a_1^{p^{m_1}},a_2^{p^{m_2}}\rg$.
  Since $\lg
  \bar{a}_1,\bar{a}_3\rg$ and $\lg\bar{a}_2,\bar{a}_3\rg$ are not metacyclic, $\lg a_1,a_2\rg$ and $\lg a_1,a_3\rg$ are not
metacyclic. By Theorem \ref{thm=d(G)=2--metacyclic}, $[a_1,a_2]^p=1$
and $[a_1,a_3]^p=1$. Moreover, $\exp(G')=p$, a contradiction.

Similar reasoning gives that $\bar{G}$ is not a group of Type (C5)--(C10) in Theorem \ref{111}.

By the above argument, $\bar{G}$ is a group of Type (C3) in
Theorem \ref{111}. That is, $\bar{G}=\bar{K}\times \bar{A}$, where
\begin{center}$\bar{K}=\lg \bar{a}_1,\bar{a}_2,\bar{a}_3\di
 \bar{a}_1^{p^{m_1+1}}=\bar{a}_2^{p^{m_2+1}}=\bar{a}_3^{p^{m_3}}=1,
  [\bar{a}_1,\bar{a}_2]=\bar{a}_1^{p^{m_1}}, [\bar{a}_1,\bar{a}_3]=\bar{a}_2^{p^{m_2}}, $

  $[\bar{a}_2,\bar{a}_3]=1\rg$, $m_1>1$ for $p=2$,

  $m_1\ge m_2\ge m_3$ and
  $\bar{A}$ is abelian such
that $\exp(\bar{A})\le p^{m_2}$.
\end{center}
Then $G'=\lg
a_1^{p^{m_1}},a_2^{p^{m_2}}\rg$.

Since $G'\not\le \lg a_2,a_3\rg$, $[a_2,a_3]=1$. Since $\lg
\bar{a}_1,\bar{a}_3\rg$ is not metacyclic, $\lg a_1,a_3\rg$
is not metacyclic. By Theorem \ref{thm=d(G)=2--metacyclic},
$[a_1,a_3]^p=1$. Let $[a_1,a_3]=a_2^{p^{m_2}}d$ where $d\in
\mho_1(G')$. Then $a_2^{p^{m_2+1}}d^p=1$. It follows that
$a_2^{p^{m_2+1}}\in\mho_2(G')$. Hence
\begin{center}
$o(a_1)>p^{m_1+1}$,
$N=\mho_1(G')=\lg a_1^{p^{m_1+1}},a_2^{p^{m_2+1}}\rg=\lg
a_1^{p^{m_1+1}}\rg$ and $a_2^{p^{m_2+1}}\in\lg a_1^{p^{m_1+2}}\rg$.
\end{center}
Since $G'\not\le \lg a_2,a_3a_1^{p^{m_1}}\rg$, $[a_2,a_3a_1^{p^{m_1}}]=1$ and hence
$[a_1^{p^{m_1}},a_2]=a_1^{p^{2m_1}}=1$. Assume that  the order of $a_1$ is
$p^{m_1+1+k}$ where $k\ge 1$. Then $m_1>k$.

Let $\bar{A}=\lg \bar{a}_4\rg\times\lg
\bar{a}_5\rg\times\dots\times\lg \bar{a}_f\rg$ and the type of
$\bar{A}$ be $(p^{m_4},p^{m_5},\dots,p^{m_f})$. For $4\le i\le f$
and $1\le j\le f$, since $G'\not\le \lg a_i,a_j\rg$, $[a_i,a_j]=1$ and hence $a_i\in Z(G)$. Assume that
$a_i^{p^{m_i}}=a_1^{sp^{m_1+1}}$. Then
$(a_ia_1^{-sp^{m_1+1-m_i}})^{p^{m_i}}=1$. Let
$b_i=a_ia_1^{-sp^{m_1+1-m_i}}$, $A=\lg b_4\rg\times\lg
b_5\rg\times\dots\times\lg b_f\rg$ and $K=\lg a_1,a_2,a_3\rg$. Then
$G=K\times A$.

Assume that $[a_1,a_2]=a_1^{p^{m_1}}a_1^{up^{m_1+1}}$. Then
$a_1^{a_2}=a_1^{1+(1+up)p^{m_1}}$. By Lemma \ref{lem=num p odd} and
Lemma \ref{lem=num p 2}, there exists a positive integer $w$ such that
$(1+(1+up)p^{m_1})^j= 1+p^{m_1}$. Replacing $a_2$ and $a_3$ with
$a_2^w$ and $a_3^w$ respectively, we have $[a_1,a_2]=a^{p^{m_1}}$.

By Lemma \ref{lem=num p odd} and Lemma \ref{lem=num p 2}, in the
multiplicative group consisting of all invertible elements of $\Bbb
Z/p^{m_1+1+k}\Bbb Z$, the order of $1+p^{m_1}$ is $p^{k+1}$. Since
$a_1^{(1+p^{m_1})^{p^{m_2+1}}}=a_1^{a_2^{p^{m_2+1}}}=a_1$, we have
$k+1\le m_2+1$. Hence $k\le m_2$.

\medskip

Case 1: $k=m_2$.

\medskip

In this case, $m_1>m_2$ and $[a_1,a_2^{p^{m_2}}]\neq 1$. It follows
that $c(\lg a_1,a_3\rg)>2$, Corollary \ref{cor=c=3 exp(G')=p} gives
that $p>2$ and $\lg a_1,a_3\rg\in \mathcal{A}_2$. If $m_3>1$, then
$\lg a_1,a_2^{p^{m_2}}a_3^p\rg$ is neither abelian nor normal in
$G$, a contradiction. Hence we have $m_3=1$. If $A\neq 1$, then,
letting $1\ne e\in A$, $\lg a_1,a_2^{p^{m_2}}e\rg$ is neither
abelian nor normal in $G$, a contradiction. Hence we have $A=1$.
Assume that $a_3^p=a_1^{vp^{m_1+1}}$. Replacing $a_3$ with
$a_3a_1^{-vp^{m_1}}$, we have $a_3^p=1$.

Assume that $[a_1,a_3]=a_2^{p^{m_2}}a_1^{wp^{m_1+1}}$. Then
$a_2^{p^{m_2+1}}a_1^{wp^{m_1+2}}=1$. Since
$$(a_2a_1^{wp^{m_1-m_2+1}})^{p^{m_2+1}}=1,$$ we may assume that
$$(a_2a_1^{wp^{m_1-m_2+1}})^{p^{m_2}}=a_2^{p^{m_2}}a_1^{wp^{m_1+1}}a_1^{xp^{m_1+m_2}}.$$
Replacing $a_2$ with $a_2a_1^{wp^{m_1-m_2+1}}a_1^{-xp^{m_1}}$, we
have $a_2^{p^{m_2+1}}=1$ and $[a_1,a_3]=a_2^{p^{m_2}}$. Hence $G$ is
a group of Type (G1) in Theorem \ref{222}.

\medskip

Case 2: $k<m_2$.

\medskip

In this case $[a_1,a_2^{p^{m_2}}]=1$. Since $[a_1,a_3,a_1]=1$, $[a_1^p,a_3]=[a_1,a_3]^p=1$. Since $G'\not\le \lg
a_2,a_3a_1^{p^{m_1-m_3+1}}\rg$,
$[a_2,a_3a_1^{p^{m_1-m_3+1}}]=1$. It follows that
$$1=[a_1^{p^{m_1-m_3+1}},a_2]=a_1^{p^{2m_1-m_3+1}}.$$ Hence
$2m_1-m_3+1\ge m_1+1+k$. That is, $m_1-m_3\ge k$. Since $G'\not\le
\lg a_1,a_2^{p^{m_2-m_3+2}}a_3^p\rg$,
$[a_1,a_2^{p^{m_2-m_3+2}}a_3^p]=1$. It follows that
$$a_1^{a_2^{p^{m_2-m_3+2}}}=a_1^{(1+p^{m_1})^{p^{m_2-m_3+2}}}=a_1.$$
By Lemma \ref{lem=num p odd} and Lemma \ref{lem=num p 2}, in the
multiplicative group consisting of all invertible elements of $\Bbb
Z/p^{m_1+1+k}\Bbb Z$, the order of $1+p^{m_1}$ is $p^{k+1}$. Hence
we have $m_2-m_3+2\ge k+1$. That is, $k\le m_2-m_3+1$.

Let $b\in A$ and the order of $b$ be $p^e$. Since $G'\not\le \lg
a_1,a_2^{p^{m_2-e+1}}b\rg$, $[a_1,a_2^{p^{m_2-e+1}}b]=1$. It
follows that
$a_1^{a_2^{p^{m_2-e+1}}}=a_1^{(1+p^{m_1})^{p^{m_2-e+1}}}=a_1$. By
Lemma \ref{lem=num p odd} and Lemma \ref{lem=num p 2}, in the
multiplicative group consisting of all invertible elements of $\Bbb
Z/p^{m_1+1+k}\Bbb Z$, the order of $1+p^{m_1}$ is $p^{k+1}$. Hence
we have $m_2-e+1\ge k+1$. That is, $e\le m_2-k$. By the
arbitrariness of $b$, $\exp(A)\le p^{m_2-k}$.

Assume that $a_3^p=a_1^{vp^{m_1+1}}$. Replacing $a_3$ with
$a_3a_1^{-vp^{m_1}}$, we have $a_3^p=1$.

Assume that $[a_1,a_3]=a_2^{p^{m_2}}a_1^{wp^{m_1+1}}$. Then
$a_2^{p^{m_2+1}}a_1^{wp^{m_1+2}}=1$. Replacing $a_2$ with
$a_2a_1^{wp^{m_1-m_2+1}}$, we have $a_2^{p^{m_2+1}}=1$ and
$[a_1,a_3]=a_2^{p^{m_2}}$. Hence $G$ is a group of Type (G2) in
Theorem \ref{222}.
 \qed

 \medskip

Summarizing, we have the following

{\bf Main Theorem.} Suppose that $G$ is a finite metahamiltonian
$p$-group. If $\exp(G')=p$, then $G$ is one of the groups listed in
Theorem \ref{111}. If $\exp(G')>p$, then $G$ is  one of the groups
listed in Theorem \ref{222}.

 \medskip

%{\bf Acknowledgments} We cordially thank Professor Joseph Brennan
%for his assistance on the language of the paper.
% Thanks go to referee for careful reading and comments. These improve the whole paper.

\end{document}